\documentclass[runningheads]{llncs}
\usepackage[T1]{fontenc}
\usepackage{graphicx}
\usepackage{xcolor}
\usepackage{comment} 
\usepackage{lipsum} 
\usepackage{amsmath}      
\usepackage{amssymb}
\usepackage{bm}           
\usepackage{algorithm}    
\usepackage{algpseudocode}
\algrenewcommand\algorithmicrequire{\textbf{Input:}}
\algrenewcommand\algorithmicensure{\textbf{Output:}}
\usepackage{bbm}
\usepackage{pgfplots}
\usepgfplotslibrary{groupplots}
\pgfplotsset{compat=1.18}
\usepackage[section]{placeins}
\usepackage{hyperref}

\usepackage{amsmath,amsfonts,bm}

\def\eqref#1{equation~\ref{#1}}

\def\1{\bm{1}}

\def\vb{{\bm{b}}}

\def\vg{{\bm{g}}}
\def\vh{{\bm{h}}}

\def\vw{{\bm{w}}}
\def\vx{{\bm{x}}}

\def\vz{{\bm{z}}}

\def\evg{{g}}
\def\evh{{h}}

\def\evx{{x}}

\def\evz{{z}}

\def\mI{{\bm{I}}}

\def\mW{{\bm{W}}}

\DeclareMathAlphabet{\mathsfit}{\encodingdefault}{\sfdefault}{m}{sl}
\SetMathAlphabet{\mathsfit}{bold}{\encodingdefault}{\sfdefault}{bx}{n}

\def\sL{{\mathbb{L}}}

\def\sN{{\mathbb{N}}}

\def\sX{{\mathbb{X}}}

\usepackage{enumerate}

\newcommand{\pga}[0]{\texttt{PGA}}
\newcommand{\ppga}[0]{\texttt{PPGA}}
\newcommand{\ppgalr}[0]{\texttt{PPGA$_{\text{LR}}$}}
\newcommand{\sw}[0]{\texttt{SimplexWalk}}

\usepackage[export]{adjustbox} 

\usepackage[outdir=./]{epstopdf}

\definecolor{lrFive}{HTML}{1f77b4}
\definecolor{lrFifty}{HTML}{ff7f0e}
\definecolor{lrFiveHundred}{HTML}{2ca02c}
\definecolor{lrFiveThousand}{HTML}{d62728}
\definecolor{lrFiftyThousand}{HTML}{9467bd}
\begin{document}
\title{Optimization over Trained Neural Networks: \\ Going Large with Gradient-Based Algorithms}
\titlerunning{Going Large with Gradient-Based Algorithms}
%
\author{Jiatai Tong\inst{1} \and
Yilin Zhu\inst{2} \and
Thiago Serra\inst{2} \and
Samuel Burer\inst{2}}
\authorrunning{J. Tong et al.}
%
\institute{Northwestern University, Evanston IL, United States\\
\email{jiataitong2026@u.northwestern.edu}
 \and
University of Iowa, Iowa City IA, United States\\
\email{\{yilin-zhu,thiago-serra,samuel-burer\}@uiowa.edu}}
\maketitle              
\begin{abstract}
%

When optimizing a nonlinear objective, one can employ a neural network as a surrogate for the nonlinear function. However, the resulting optimization model can be time-consuming to solve globally with exact methods. As a result, local search that exploits the neural-network structure has been employed to find good solutions within a reasonable time limit. For such methods, a lower per-iteration cost is advantageous when solving larger models. The contribution of this paper is two-fold. First, we propose a gradient-based algorithm with lower per-iteration cost than existing methods. Second, we further adapt this algorithm to exploit the piecewise-linear structure of neural networks that use Rectified Linear Units~(ReLUs). In line with prior research, our methods become competitive with---and then dominant over---other local search methods as the optimization models become larger.

\keywords{Constraining learning \and Gradient ascent \and Linear regions \and Piecewise-linear functions \and Rectified linear units.}
\end{abstract}
\section{Introduction}

In the field of mathematical programming, 
piecewise-linear functions play an important role in modeling nonlinear functions ~\cite{mccormick1976under,mangasarian2005under,vielma2010framework,misener2010multi,huchette2023pwl}. 
In deep learning, a popular model that provides a piecewise-linear approximation of a nonlinear function is the neural network with the ReLU activation function~\cite{arora2018understanding,huchette2023survey}.
Researchers have long known that 
some neural network architectures are universal function approximators~\cite{cybenko1989approximation,funahashi1989approximate,hornik1989approximator}, 
and in particular, this is also true of the ReLU activation function if the architecture is sufficiently wide~\cite{yarotsky2017relu} or deep~\cite{hanin2017approximating}. 
When neural-network approximations are used as surrogates for solving nonlinear optimization problems, algorithms that exploit the piecewise-linear structure of the neural networks are of particular interest.
In this paper, 
we propose gradient-based algorithms for this setting. 

Optimizing a piecewise-linear function over a polyhedron can be modeled using a Mixed-Integer Linear Programming~(MILP) formulation. In the specific case of ReLU networks,
existing MILP formulations either have a weak linear relaxation due to big M coefficients~\cite{fischetti2018constraints} or become prohibitively large when using a disjunctive formulation~\cite{huchette2018phd}. Researchers have found success by improving the big $M$ coefficients~\cite{cheng2017resilience,grimstad2019surrogate,liu2021algorithms,badilla2023tradeoff,zhao2024horizon,hojny2024message,sosnin2024control}, 
strengthening formulations using valid inequalities~\cite{anderson2020strong}, 
and using a hybrid of both formulations~\cite{tsay2021partition}.
Other improvements include reformulation~\cite{huchette2018phd,schweidtmann2019global,liu2025doc}, 
parameter rescaling~\cite{plate2025scaling}, 
pruning the search space by inference~\cite{tjeng2019evaluating,xiao2019training,botoeva2020efficient,serra2020empirical}, and 
working with sparser neural networks~\cite{say2017planning,xiao2019training,cacciola2024structured,pham2025surrogate}. 
Those improvements can also help using MILP formulations in neural networks for verification~\cite{anderson2020strong,cheng2017resilience,rossig2021verification,strong2021global},  compression~\cite{serra2021compression,serra2020lossless}, evaluating output variation~\cite{cai2023pruning,kumar2019equivalent,liu2020monotonic,serra2020empirical,serra2018bounding}, and counterfactual explanations~\cite{kanamori2021counterfactual,tsiourvas2024counterfactual}.

\begin{figure}
    \centering
    \includegraphics[width=0.96\linewidth]{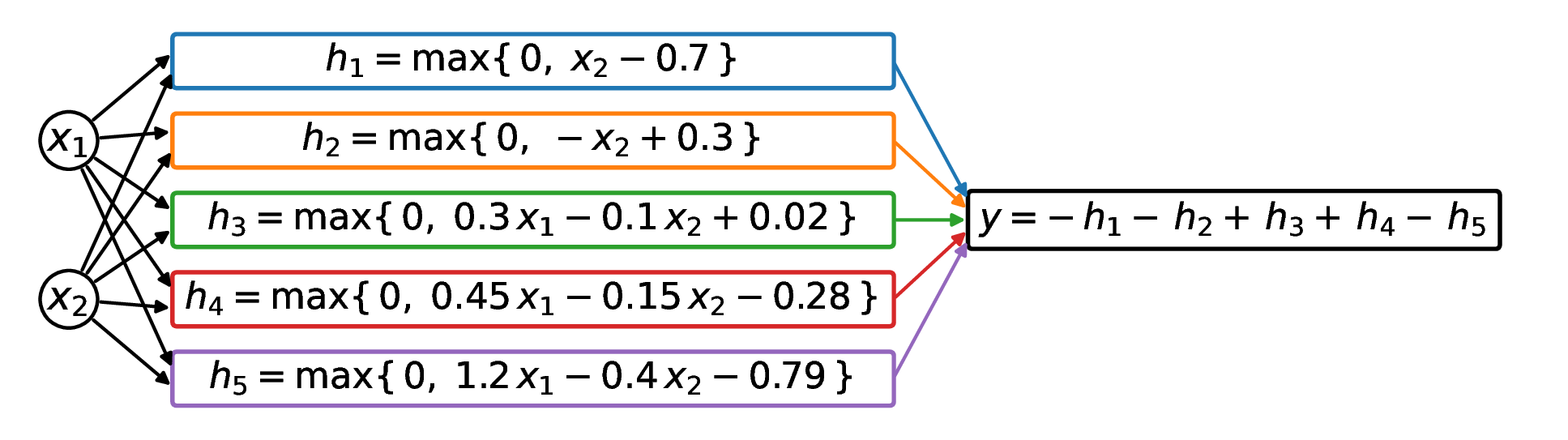}
    \includegraphics[scale=0.48]{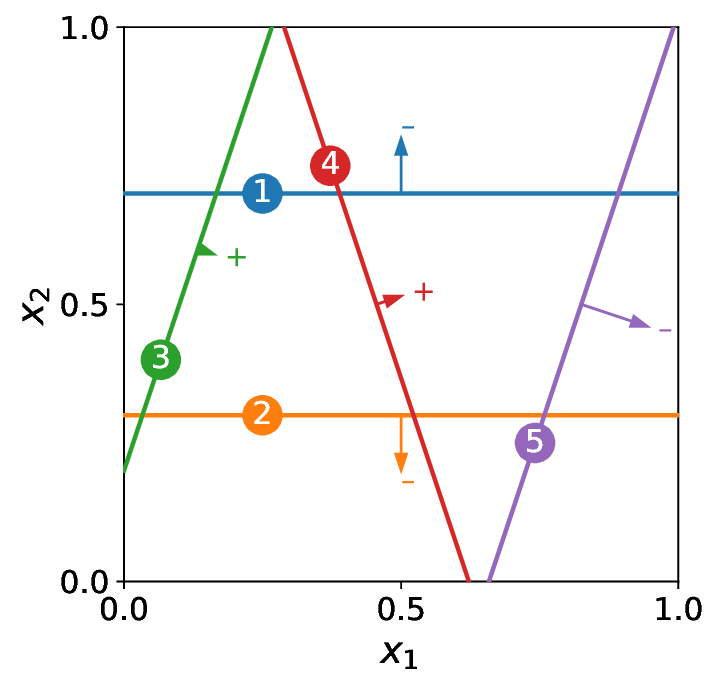}
    \includegraphics[scale=0.48]{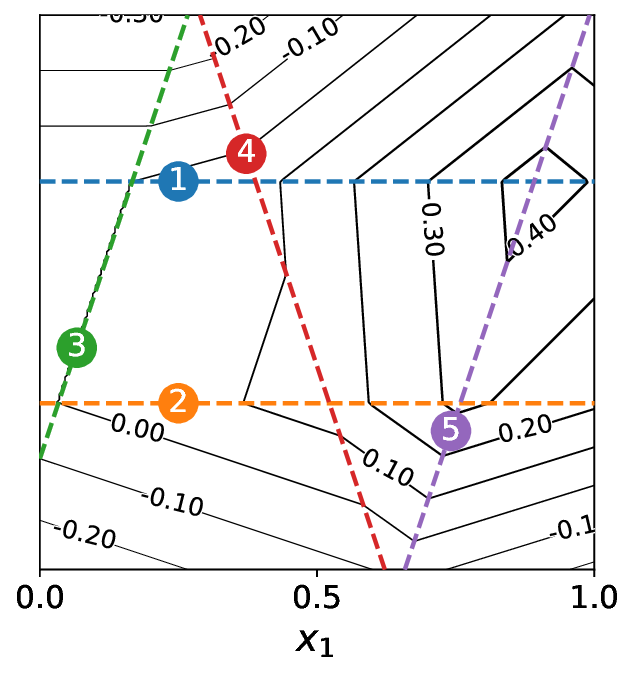}
    \caption{
    \textbf{Top:} Visual description of neural network used as example. 
    \textbf{Bottom left:} Lines partitioning the space based on what inputs produce a positive output for each neuron, with the arrow pointing to the positive side, the length of the arrow proportional to the magnitude of the parameters, and the arrow label denoting the influence on $y$.
    \textbf{Bottom right:} Contour plots of $y$ over the lines associated with neural activations.
    }
    \label{fig:pwa}
\end{figure}

There are also local search methods, which do not solve MILPs exactly, but instead are designed to find good solutions in limited time. Indeed, both Perakis \& Tsiourvas~\cite{perakis2022optimizing} and Tong et al. \cite{tong2024optimization} are closely related to this paper for taking a geometric view of the input space of the ReLU network, focusing on the linear pieces of the function approximation. Known as \emph{linear regions} in machine learning, each piece corresponds to a polyhedron associated with a distinct set of active neurons. Within a linear region, changes to the input have a linear impact on the output. In Figure~\ref{fig:pwa}, we illustrate these concepts 
on a neural network having inputs $x_1$ and $x_2$, 
neurons with outputs $h_1$ to $h_5$, 
and output $y$.

\begin{figure}
    \centering
    \includegraphics[scale=0.49,valign=b]{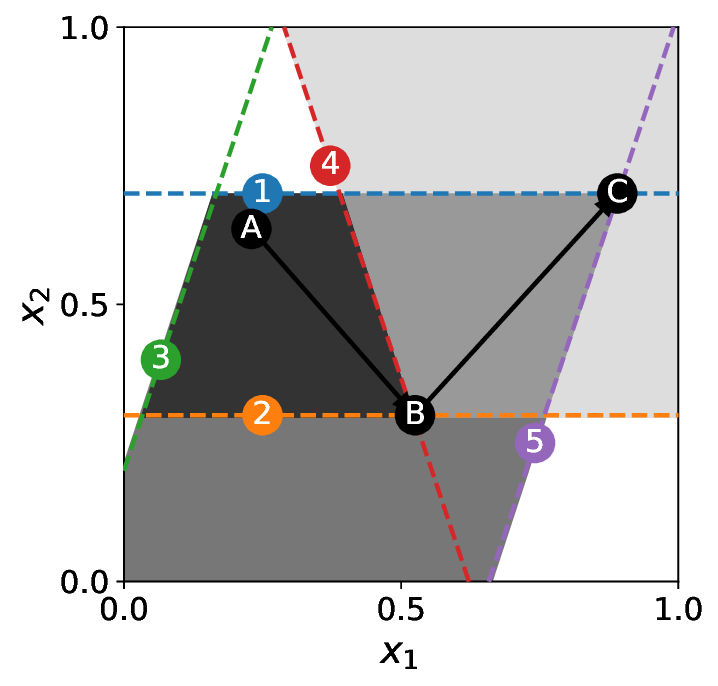}
    \includegraphics[scale=0.49,valign=b]{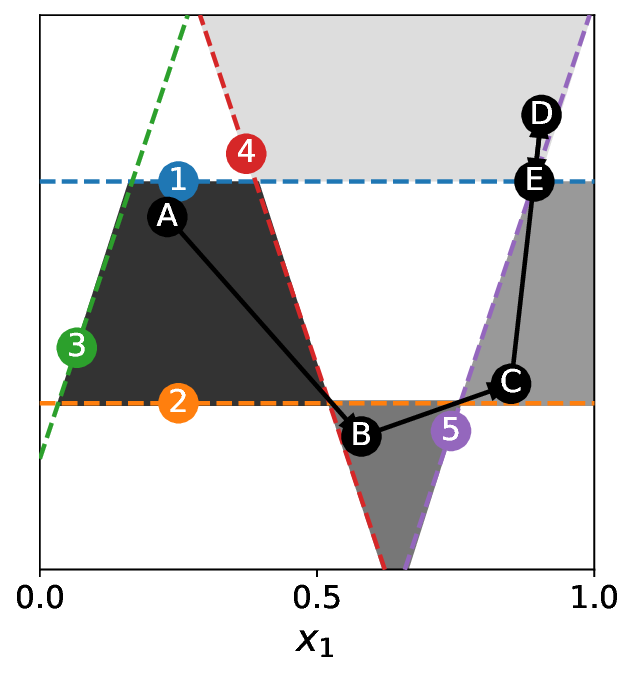}
    \caption{\textbf{Left:} Solutions found by MILP Walk from the initial solution $A = (0.23, 0.636)$ until convergence. \textbf{Right:} Solutions found by LP Walk from $A$ until convergence.}
    \label{fig:other_walks}
\end{figure}

We summarize the works \cite{perakis2022optimizing,tong2024optimization} just mentioned, taking the liberty to name them {\em MILP Walk\/} and {\em LP Walk\/}, respectively, in order to draw parallels between these methods and our methods introduced in Section \ref{sec:algorithms}:
\begin{itemize}
    \item \textbf{MILP Walk:} Perakis \& Tsiourvas~\cite{perakis2022optimizing} solve a sequence of restricted MILP models. 
    Each MILP model finds the best solution across all linear regions containing the current solution. 
    If a better solution is found, the same process is repeated from the new solution.  
    In Figure~\ref{fig:other_walks} Left, solution $A$ lies only in the linear region in darkest gray, in which the best solution is $B$. 
    In turn, solution $B$ lies in the four linear regions with the three darker tones of gray, 
    where the best solution is $C$. 
    Finally, solution $C$ lies in the four linear regions with the two lighter tones of gray, 
    where the best solution is $C$ again. Once no improvement is found, the algorithm stops.
    \item \textbf{LP Walk:} Tong et al.~\cite{tong2024optimization} solve a sequence of LP models. 
    Each LP model finds the best solution in a linear region containing the current solution. 
    If a better solution is found, they repeat the process by moving slightly past the new solution along the line from the last solution. Moving slightly past avoids using a solution lying in multiple linear regions. 
    In Figure~\ref{fig:other_walks} Right, the linear region of solution $A$ is in darker gray, and its best solution is along the line from $A$ to $B$. 
    In turn, the linear region of solution $B$ is in slightly lighter gray, and its the best solution is along the line from $B$ to $C$. 
    From $C$ we find solution $E$ and move towards $D$. 
    From $D$ we find solution $E$ again. At this point, the algorithm stops.
\end{itemize}
For a broader discussion about linear regions and their nexus with mathematical optimization, 
we recommend the survey by Huchette et al.~\cite{huchette2023survey}.

Which algorithm, MILP Walk or LP Walk, is best suited for a particular instance often depends on the size of that instance. Starting from the same solution, it is easy to see that MILP Walk will move next to a solution that is at least as good as the one found by LP Walk. On the other hand, each iteration of MILP Walk solves the same MILP model used to optimize over the entire neural network, albeit restricted to the neighborhood of the current solution. Hence,
the per-iteration cost of MILP Walk is higher than the the per-iteration cost of LP Walk. Consequently, LP Walk can perform more iterations in the same amount of time. Indeed, an empirical comparison of both methods shows that LP Walk performs better for neural networks with more inputs, layers, and neurons~\cite{tong2024optimization}, all of which imply a larger number of linear regions~\cite{serra2018bounding}. Thus, LP Walk conducts a style of search that is more akin to sampling than to enumeration~\cite{serra2020empirical}. 

\begin{figure}
    \centering
    \includegraphics[scale=0.49,valign=b]{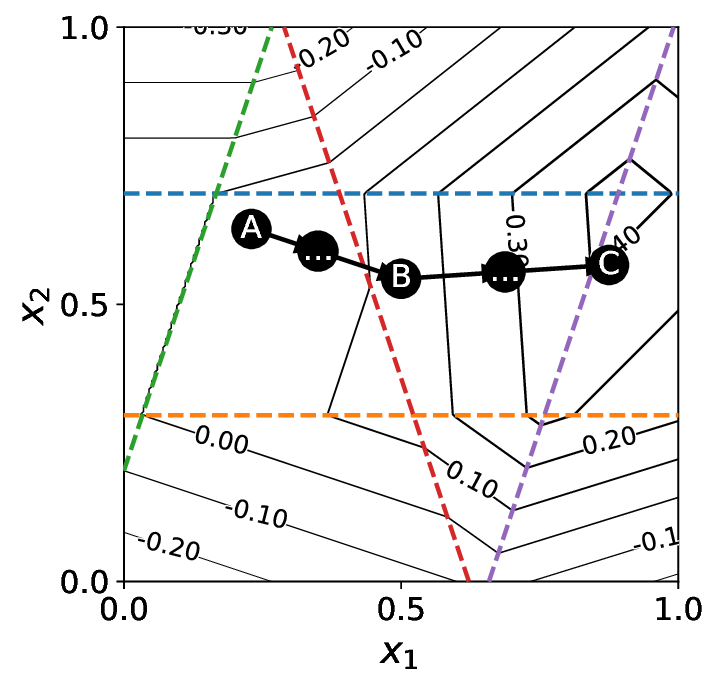}
    \includegraphics[scale=0.49,valign=b]{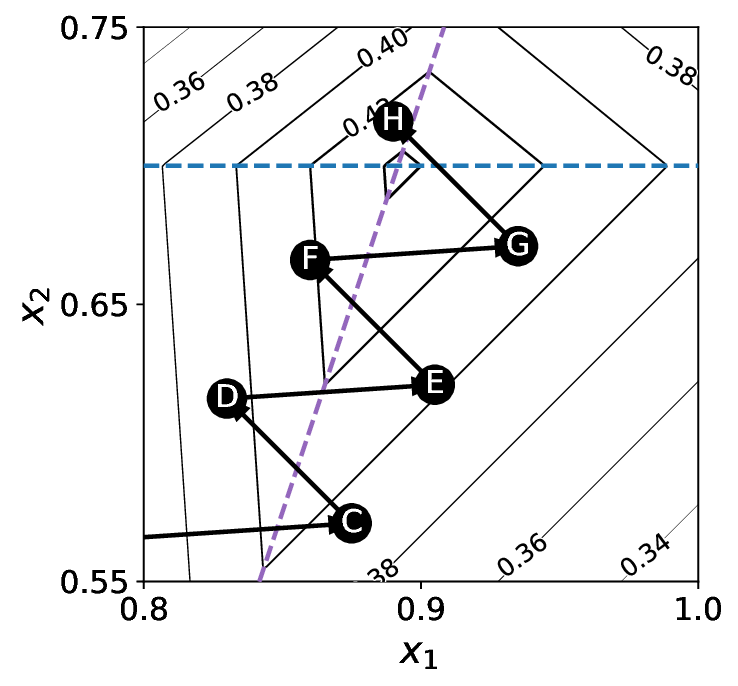}
    \caption{\textbf{Left:} First solutions found by Gradient Walk from the initial solution $A = (0.23, 0.636)$. \textbf{Right:} Next solutions found over a narrower region of the input space.}
    \label{fig:our_walk}
\end{figure}

Of course, solving an LP model for each linear region, as in LP Walk, may eventually become too costly in ever larger neural networks. Hence, in this paper, we propose a new local search approach with an even smaller per-iteration cost:
\begin{itemize}
    \item \textbf{Gradient Walk:} We compute a sequence of gradient steps. 
    Each step may find a better solution within the current linear region, or a solution in another linear region that is better or worse. We keep track of the best solution found. 
    If the improvement is too small for a predefined number of steps, 
    we restart from a perturbation of the best solution found thus far. 
    In Figure~\ref{fig:our_walk} Left, 
    we move from $A$ to $B$ through nine intermediary steps in the same linear region as $A$, all in the same direction, 
    and then from $B$ to $C$ with five intermediary steps. 
    Note that the steps are orthogonal to the contour plots within the linear regions. 
    In Figure~\ref{fig:our_walk} Right, 
    we continue from $C$ until $H$ by zig-zagging among linear regions 
    and improving in all steps except $G$. 
    The  steps following H would not find a much better solution.
    Hence, the algorithm stops.
\end{itemize}

In what follows, 
we define our problem of interest and its conventional MILP model in Section~\ref{sec:notation}. 
Then we present an algorithm for Gradient Walk as well as a variant that further exploits knowledge of linear regions in Section~\ref{sec:algorithms}. We evaluate those algorithms in Section~\ref{sec:exp}. Conclusions are given in Section \ref{sec:conclusion}.

\section{Preliminaries}\label{sec:notation}

We optimize the function $f : \mathbb{R}^{n_0} \rightarrow \mathbb{R}$ associated with a neural network over a polytope $\mathbb{X}\subset \mathbb{R}^{n_0}$:
\begin{align}
    \mathrm{maximize}_\vx ~~ &f(\vx) \label{eq:opt_obj}\\
    \mathrm{s.t.} ~~ &\vx\in\mathbb{X} \label{eq:opt_constr}
\end{align}
i.e., we want an input $\vx = [\evx_1 ~ \evx_2 ~ \dots ~ \evx_{n_0}]^\top \in \mathbb{X}$ maximizing the prediction $f(\vx)$. In the experiments of Section \ref{sec:exp}, $\mathbb{X}$ equals a box, but our theoretical development requires only that $\mathbb{X}$ be polyhedral.

Let the neural network have $L$ hidden layers, each hidden layer $l \in \mathbb{L} := \{1, \dots,L\}$ having preactivation values $\vg^l =[\evg_1^l ~ \evg_2^l \dots \evg_{n_l}^l]^\top$ and outputs $\vh^l = [\evh_1^l ~ \evh_2^l \dots \evh_{n_l}^l]^\top$ from neurons indexed by $i \in \sN_l = \{1, 2, \ldots, n_l\}$. The output layer $L+1$ has a single preactivation value $\evg_1^{L+1}$ as the network output.
Let $\mW^l$ and $\vb^l$ denote the weight matrix and bias vector associated with the $l$-th layer, 
for which we assume that all elements are within $[-1,1]$. The preactivation value $\evg^l_i$ of neuron $i$ in layer $l$ is given by $\evg_i^l = \mW_{i}^l \vh^{l-1} + \vb_i^l$, and the output $\evh^l_i$ of neuron $i$ in hidden layer $l$ follows by the ReLU activation $\evh^l_i=\max\{0, \evg^l_i\}$. In this setting, $\vx = \vh^0$ and $f(\vx) = \evg^{L+1}_{1}$.

For an input $\vx \in \sX$, let $\vz^l(\vx) = [\evz_1^l ~ \evz_2^l \dots \evz_{n_l}^l]^\top$ be the \emph{layer activation pattern} produced in layer $l$ of the neural network when given $\vx$ as input, where:
\begin{equation}\label{eq:activation_pattern_cases}
z_i^{\,l} =
\begin{cases}
1, & \text{if } h_i^{\,l} = g_i^{\,l} \geq 0,\\
0, & \text{if } g_i^{\,l} \leq 0.
\end{cases}
\end{equation}
A neuron is \emph{binding} if $\evg^l_i = \evh^l_i = 0$. 
In this case, $\evz^l_i$ can be either $0$ or $1$. 

Let $\vz(\vx) = \{ \vz^1(\vx), \vz^2(\vx), \ldots, \vz^L(\vx)\}$ be the \emph{activation pattern} produced across all layers of the neural network when given $\vx$ as input. 
A \emph{linear region} $\mathcal{R}_{\vz'} \in \sX$ is a set where every input $\vx$ has the same activation pattern $\vz'$ as all the other inputs, 
i.e., $\vz(\vx) = \vz' ~ \forall \vx \in \mathcal{R}_{\vz'}$. 
The output of $f$ varies linearly within each linear region. 
Inputs lie in multiple linear regions if a neuron is binding.

We can formulate the optimization problem as the following MILP model:
\begin{align}
    \mathrm{maximize} ~~ & f(\vx) = g_1^{L+1} \\
    \mathrm{s.t.} ~~ & \vh^0 = \vx, \qquad \vx \in \sX \\
    & \mW_i^l \vh^{l-1} + \vb_i^l = \vg_i^l, & \forall l \in \sL \cup \{ L+1 \}, i \in \sN_l \label{eq:mip_unit_begin} \\
    & (\vz_i^l = 1) \rightarrow (\vh_i^l = \vg_i^l), & \forall l \in \sL, i \in \sN_l \label{eq:first_indicator} \\  
    & (\vz_i^l = 0) \rightarrow (\vg_i^l \leq 0 \wedge \vh_i^l = 0), & \forall l \in \sL, i \in \sN_l \label{eq:last_indicator} \\
    & \vh_i^l \geq 0 \label{eq:lp_unit_end}, \qquad \vz^l_i \in \{0,1\}, & \forall l \in \sL, i \in \sN_l 
\end{align}
The indicator constraints (\ref{eq:first_indicator})--(\ref{eq:last_indicator}) can be modeled with big M constraints~\cite{bonami2015indicator}. 

We can formulate an LP model optimizing over a single linear region $\mathcal{R}_{\vz'}$ by fixing $\vz = \vz'$. For example, this is what LP Walk does in each iteration.

\section{Our Proposed Algorithms}\label{sec:algorithms}

We propose two gradient-based algorithms for solving (\ref{eq:opt_obj})--(\ref{eq:opt_constr}). We assume throughout that the polyhedral domain $\mathbb{X}$ is bounded, i.e., $\mathbb{X}$ is a polytope. Indeed, neural networks are often used within a confined domain. There is a growing body of work on preventing a neural network from extrapolating outside the nominal or implied domain defined by a training set~\cite{chi2024trust,tsiourvas2024counterfactual,zhu2024validity}. 
Moreover, 
the possibility of embedding a neural network as part of an MILP model depends on the input set being bounded~\cite{serra2018bounding}, 
which otherwise would also prevent us from benchmarking against existing algorithms~\cite{perakis2022optimizing,tong2024optimization}.

Our first algorithm, the Perturbed Projected Gradient Ascent~(\ppga{}) algorithm, is described in Section~\ref{sec:PPGA}. \ppga{} makes use of both the piecewise-constant nature of $\nabla f$ and projection onto $\mathbb{X}$. Our second algorithm, \ppga{} with Linear Region Valve~(\ppgalr{}), is described in Section \ref{sec:PPGALR}. \ppgalr{} enhances \ppga{} by further exploiting the structure of the linear regions, being particularly beneficial if the number of linear regions is large.

\subsection{Perturbed Projected Gradient Ascent (Algorithm \ref{alg:ppga})} \label{sec:PPGA}

We first mention the standard approach, called the Projected Gradient Ascent~(\pga{}) algorithm, of projecting the gradient steps over the feasible set to produce iterates of the form 
\begin{align}
\vx^{t+1} = P_\mathbb{X}(\vx^t + \gamma\nabla f(\vx^t)),
\end{align}
where $\gamma$ is the learning rate and $L_2$ projection~\cite{rosen1960GradientProjection,polyak1966ProjectionRroof} 
solves the convex quadratic program (QP)
\begin{align}
    P_\mathbb{X}(\dot{\vx}) = \arg \min_\vx ~~ &\|\vx - \dot{\vx}\|^2 \label{eq:proj_obj}\\
    \mathrm{s.t.} ~~ &\vx\in\mathbb{X}, \label{eq:proj_constr}
\end{align}
which takes only $O(n_0)$ time when $\sX$ is a box.

Based on \pga{}, we now introduce our method \ppga{}. For any point $\vx' \in \mathcal{R}_{\vz(\vx)}$, i.e., any point $\vx'$ in the same linear region as $\vx$, 
we can calculate $f(\vx')$ with the affine transformation 
\begin{align}
f|_{\mathcal{R}_\vz(\vx)}(\vx')
= T(\vz(\vx)) \vx' + t(\vz(\vx)),
\end{align}
where 
\begin{align}
T(\vz(\vx)) = \vw^{L+1} \left( \prod_{l = 1}^L (\vz^l(\vx)\mI)\mW^l \right)
\end{align}
and $t(\vz(\vx))$ is constant~\cite{huchette2023survey}. 
Hence, 
it follows that 
\begin{align}
    \nabla f(\vx) = \vw^{L+1}\left(\prod_{l = 1}^L (\vz^l\mI)\mW^l\right) \label{eq:grad}
\end{align}
for any point $\vx$ at the interior of its linear region, i.e., not binding for any neuron. 
Hence, $f$ has a piecewise affine landscape,
which may contain local maxima, saddle points, and local minima. 

\begin{algorithm}[b!]
\caption{Perturbed Projected Gradient Ascent
}
\label{alg:ppga}
\begin{algorithmic}[1]
  \Require Model $f$ with input size $n_0$, learning rate $\gamma$, time limit $T$,  
           restart noise coefficient $\Xi$, 
           error threshold $\epsilon$,  
           tolerance window $k$; input domain $\mathbb{X}$ 
  \Ensure Best solution $\vx^{*}$ and objective value $f(\vx^*)$
  \State $\delta \gets \dfrac{\Xi}{\sqrt{n_0}}$
  \State Sample initial $\vx \sim \mathrm{Uniform}(\mathbb{X})$  
  \State $\vx^* \gets \vx$ \Comment{Initialized best solution so far}
  \State $\vx' \gets \vx$ \Comment{Initialized best solution since reset}
  \State $r \gets 0$ \Comment{Initialized reset counter}
  \While{Time $< T$}
    \State $\vx \gets P_\mathbb{X}\left( \vx + \gamma \nabla f(\vx) \right)$ \label{lin:grad_step} \Comment{Gradient step; replaced with Algorithm~\ref{alg:PPGALR} in \ppgalr{}}
    \If{$f(\vx) > f(\vx')$} \Comment{If found best solution since reset}
      \State $\Delta \gets f(\vx) - f(\vx')$ \Comment{Calculate local improvement before update}
      \State $\vx' \gets \vx$ \Comment{Update best solution since reset}
      \If{$f(\vx) > f(\vx^*)$} \Comment{If found best overall solution}
        \State $\vx^* \gets \vx$ \Comment{Update best overall solution}
      \EndIf
      \If{$\Delta < f(\vx) \cdot \epsilon $} \Comment{If local improvement is below threshold}
          \State $r \gets r+1$ \Comment{Increment reset counter}
          \If{$r = k$} \Comment{If reached reset trigger}
            \State $\vx \gets P_\mathbb{X}\bigl(\vx^* + \xi\bigr),\quad \xi \sim \mathcal{N}(0,\delta)$ \Comment{Reset event}
            \State $\vx' \gets \vx$
            \State $r \gets 0$
          \EndIf   
          \Else  \Comment{If the improvement is above the threshold}
            \If{$f(\vx) = f(\vx^*)$}    \Comment{If  current solution matches best overall solution}
                \State $r \gets 0$  \Comment{Reset the counter}
            \EndIf
      \EndIf
    \EndIf
  \EndWhile
  \State \Return $\vx^{*},\;f(\vx^*)$
\end{algorithmic}
\end{algorithm} 

For smooth maximization problems, 
one may generally escape from (interior feasible) saddle points and local minima by introducing a perturbation when $\nabla f(\vx)$ is small enough \cite{jin2017saddle1,jin2018saddle2,guo2020saddle3,vahedi2024spgd}. 
Our case is quite different, however, since all differentiable regions have constant gradients  
and because saddle points and local minima occur at the boundary of those regions, where the function is nondifferentiable. Hence, the size of $\nabla f(x)$ is not a reliable measure of local optimality in our case.

We hence use a different mechanism to measure progress. If we accumulate $k$ improvements that are relatively small in comparison to the current objective value $f(x^t)$, while not producing a better overall solution, then we continue the next iteration from a perturbation $\xi \sim \mathcal{N}(0, \frac{\Xi}{\sqrt{n_0}})$ around the current best solution $x^*$.
If that happens at step $t$, then the next update is 
\begin{align}
x^{t+1} = P_{\sX}(x^* + \xi + \gamma \nabla f(x^* + \xi))
\end{align}
The resultant algorithm, our Perturbed Projected Gradient Ascent~(\ppga{}) algorithm, is described as Algorithm~\ref{alg:ppga}.

\subsection{\ppga{} With Linear Region Valve (Algorithms \ref{alg:ppga} and \ref{alg:PPGALR})}\label{sec:PPGALR}

When the ReLU network gets deeper, the gradient calculated with Equation~(\ref{eq:grad}) may either explode or diminish, 
which makes it important to find an appropriate learning rate. 
As an alternative to calibrate the learning rate, 
we propose using linear region information for making local decisions about the size of the gradient step. 
Because we assume in this paper that all weights are within $[-1, 1]$, we expect the gradient to diminish when the network gets deeper. Notably, a similar approach can also be applied to resolve gradient explosion.

Suppose that we are in the linear region $\mathcal{R}_\vz$ with activation pattern $\vz$. We can calculate how far we may move to reach the next linear region in the direction of the gradient by a ratio test:
\begin{align}
    u = \operatorname*{min}_{i\in\mathcal{I}} \bigl(-\frac{g_i}{\Delta g_i (\vx)}\bigr) \label{eq : ratio test}
\end{align}
where $\mathcal{I} = \left\{i \mid i \in \sN_l, l \in \sL, -\frac{g_i}{\Delta g_i (\vx)} \geq 0 \right\}$ is a subset of all neurons, and 
\begin{align}
    \Delta g_i(\vx) = g_i(\vx + \nabla f(\vx)) - g_i(\vx). \label{eq:empirical delta g}
\end{align}
Here, $u$ is the relative step size to the next linear region, while the actual step size is $u\cdot \|\nabla f(\vx^t)\|$. We can use $u$ as an estimate for the size of a linear region $\mathcal{R}_{\vz'}$ near $\mathcal{R}_\vz$, i.e., $\|\vz - \vz'\|_1 \leq \zeta$ with a relatively small $\zeta$. 
Given also the relative step size $\gamma$, 
then we estimate the gradient step to stretch over 
\begin{align}
    v = \left\lceil \frac{\gamma}{u} \right\rceil
\end{align}
linear regions around $\mathcal{R}_\vz$. 

Let $V > 1$ be a predetermined valve value, 
which corresponds to the number of linear regions that we would like to stretch over at each gradient step. 
If $\gamma$ is such that $v < V$ at the current linear region, 
then we use a scale factor $c$ over the magnitude of the gradient:
\begin{align}
\vx^{t+1}
= 
\begin{cases}
P_\mathbb{X}\left(\vx^t 
+ c\, \dfrac{\nabla f(\vx^t)}{\|\nabla f(\vx^t)\|}\right), 
& \text{if } \gamma \le V \cdot u, \\[1.2ex]
P_\mathbb{X}\left(\vx^t + \gamma\, \nabla f(\vx^t)\right),
& \text{otherwise}.
\end{cases} \label{eq: valve}
\end{align}

The adaptive gradient step replaces Line~\ref{lin:grad_step} in Algorithm~\ref{alg:ppga} with Algorithm \ref{alg:PPGALR}.

\begin{algorithm}[h]
    \caption{Adaptive Gradient Step with Linear Region Valve}
    \begin{algorithmic}[1]
    \Require Current solution $\vx$, learning rate $\gamma$,  valve value $V$, scale factor $c$
    \Ensure New iterate $\vx'$ \Comment{Replaces iterate calculated in Line~\ref{lin:grad_step} of Algorithm~\ref{alg:ppga}}
    \State $\nabla f(\vx) \gets \vw^{L+1}\left(\prod\limits_{l = 1}^L (\vz^l\mI)\mW^l\right)$
    \State $g \gets f(\vx)$
    \State $g' \gets f\bigl(\vx + \nabla f(\vx)\bigr)$
    \State $\Delta g \gets g' - g$
    \State $\rho \gets - g / \Delta g$
    \State Mask all negative entries in $\rho$ with $+\infty$
    \State $u \gets \min(\rho)$
    \If{$V\cdot u \ge \gamma$}
        \State $\vx' \gets P_\mathbb{X}\left(\vx + c\cdot\dfrac{\nabla f(\vx)}{\bigl\lVert  \nabla f(\vx) \bigr\rVert}\right)$ 
    \Else
        \State $\vx' \gets P_\mathbb{X}\bigl(\vx + \gamma \cdot \nabla f(\vx)\bigr)$
    \EndIf
    \State \Return $\vx'$
    \end{algorithmic}
    \label{alg:PPGALR}
\end{algorithm}

In our implementation, we simply chose $V = \frac{1}{\|\nabla f(\vx^t)\|}$ as an adaptive valve value, so that a small gradient will more likely trigger the rescaled gradient update. We also set $c = u$, so that we force the step size stretching more than $V$ linear regions, since each linear region is estimated to have length $u \cdot \|\nabla f(\vx^t)\| = \frac{u}{V}.$ The main reason using adaptive hyperparameters is that applying grid search to extra hyperparameters is very costly, and we want to be fair to those algorithms with fewer hyperparameters, such as \ppga{} and \cite{tong2024optimization}. 

We denote this variant of \ppga{} using Algorithm~\ref{alg:PPGALR} and adaptive $(V,c)$ design as \ppgalr{}.

\section{Numerical Experiments} \label{sec:exp}

We devised numerical experiments to evaluate algorithms \ppga{} and \ppgalr{} 
on standard benchmarks and compare them with other methods. In the following subsections, we define the concept of a basic experiment, then describe our method for generating multiple experiments, and finally detail the optimization results. 
All numerical experiments were implemented in Python 3.10.8 using Gurobi 11.0 and evaluated on a single Xeon E5-2680v4 core running at 2.4 GHz with 16 GB of memory under the CentOS Linux operating system. The source code is publicly shared at \href{https://github.com/yillzhu/nn_opt}{https://github.com/yillzhu/nn\_opt} .

\subsection{Definition of an Experiment}

In our study, an {\em experiment\/} $\mathcal{P}$ refers to a complete specification of five design options that determine the structure of a neural network and the algorithm used:

\begin{itemize}

    \item {\em Input size\/}, i.e., the input dimension $n_0$ of the neural network.

    \item {\em Depth\/}, i.e., the number of hidden layers $d$ of the neural network.

    \item {\em Width\/}, i.e., the number of neurons $m$ in each hidden layer.

    \item {\em Seed\/}, i.e., the seed $s$ used to instantiate the network parameters.

    \item {\em Algorithm\/}, i.e., the local search method $\mathcal{M}$ used for optimization.
    
\end{itemize}
Given a particular combination of these specifications, the experiment proceeds as follows.  
With seed $s$ fixed, we generate a neural network having input dimension $n_0$, depth $d$, and width $m$. For each optimization algorithm $\mathcal{M}$, we then conduct a grid search to determine the optimal hyperparameters: shaking noise $\sigma$, error threshold $\epsilon$, and tolerance window $k$. Once these parameters are fixed, the algorithm $\mathcal{M}$ is used to optimize the network within a specified time limit.

In line with prior work~\cite{perakis2022optimizing,tong2024optimization} and to properly benchmark with it, 
we optimize over neural networks with their weights as defined at initialization. 
Hence, instead of optimizing over neural networks approximating specific functions based on their training, 
we work with neural networks representing distinct and random functions. 
We believe that this makes the results more representative.

\subsection{Generating Multiple Experiments} \label{sec:exp:generate_exp}

In addition to \ppga{} and \ppgalr{}, 
we use the \pga{} algorithm from Section~\ref{sec:PPGA} as a baseline 
and benchmark against the algorithm proposed for LP Walk in~\cite{tong2024optimization}, which we denote as \sw{}. 
\sw{} has shown better scalability than solving directly with Gurobi~\cite{gurobi} or with MILP Walk~\cite{perakis2022optimizing}. 
Hence, 
\[
\mathcal{M} \in \left\{  \ppgalr{}, \ppga{}, \pga{}, \sw \right\}.
\]
We use input size $n_0 \in \{10, 100,1000\}$, depth $d \in \{2, 4, 6, 8\}$, and width $m \in \{100, 1000, 10000\}$, resulting in $3 \times 3 \times 4 = 36$ distinct network configurations.

For gradient-based algorithms (\ppgalr{}, \ppga{}, \pga{}), performance is very susceptible to hyperparameters. Therefore, we use a preprocessing phase to search for better hyperparameters for each gradient-based algorithm through grid search and voting. For each network configuration $(n_0, d, m)$, we generate five random networks using seeds $s \in \{5,6,7,8,9\}$. These networks, denoted as {\em grid search instances}, are used for choosing hyperparameters over a grid with
\[
\gamma \in \{0.001, 0.01, 0.1, 1, 5\}, \quad
\sigma \in \{0.2, 2, 20\}, \quad
k \in \{100, 500, 1000\}.
\]
For each grid search instance, we evaluate all $5 \times 3 \times 3 = 45$ parameter combinations, running the optimization for 5 minutes per grid point. The five combinations achieving the best objective values are recorded for each instance.  
After processing all five grid search instances, we obtain five sets of top-performing parameter combinations. The most frequently occurring combination across these sets is selected as the final grid search result for that algorithm and network configuration. This process runs independently for each gradient-based algorithm.

Once the hyperparameters are calibrated, we generate 20 additional networks for each $(n_0, d, m)$ combination using seeds $s \in \{10, 11, \ldots, 29\}$. These are referred to as {\em optimization instances}. Each algorithm $\mathcal{M}$ is then applied to these instances, and the gradient-based algorithms use the hyperparameters obtained from the grid search. Every run is executed for 7200 seconds, during which we record, at every second, the best objective value and the number of iterations.

\subsection{Computational Complexity of Walk Steps}
The complexity of the problem can be affected by all three setup options $(n_0, d, m)$. For gradient-based algorithms (\ppgalr{}, \ppga{}, \pga{}), the asymptotic complexity of both calculating the gradient and making a prediction is
$$O(dm^2 + n_0m).$$
which can be deduced from a sequence of vector-matrix multiplications. Through Equation~(\ref{eq:empirical delta g}),  each step of the linear region algorithm (\ppgalr{}) may cost twice as much as \pga{} and \ppga{} steps, 
which we may consider still acceptable. 
On the other hand, \sw{} carries out the solution of one LP model per step. 
Notably, the computational cost of a \sw{} step, if the LP model is solved using the simplex algorithm, is exponential on $n_0$ in worst-case.

\subsection{Results} \label{sec:result}

We use the Dolan–More performance profile~\cite{MR1875515,MR1413376} to compare the performance of different algorithms on multiple problems. A performance profile shows, for each algorithm, the fraction of test instances on which it performs within a given factor of the best observed result. The horizontal axis represents the performance factor $\tau \geq 1$ in log scale, and the vertical axis $\rho_{\mathcal{M}}(\tau)\in [0,1]$ represents the fraction of experiments for which the algorithm $\mathcal{M}$ attains a performance ratio of at most $\tau$. The vertical intercept $\rho_{\mathcal{M}}(1)$ indicates how often an algorithm achieves the best result among all competitors, and a curve approaching $\rho = 1$ more rapidly reflects an algorithm whose performance is consistently close to the best algorithm across all experiments.

\begin{figure}[h!]
    \centering
    \includegraphics[width=1.0\linewidth]{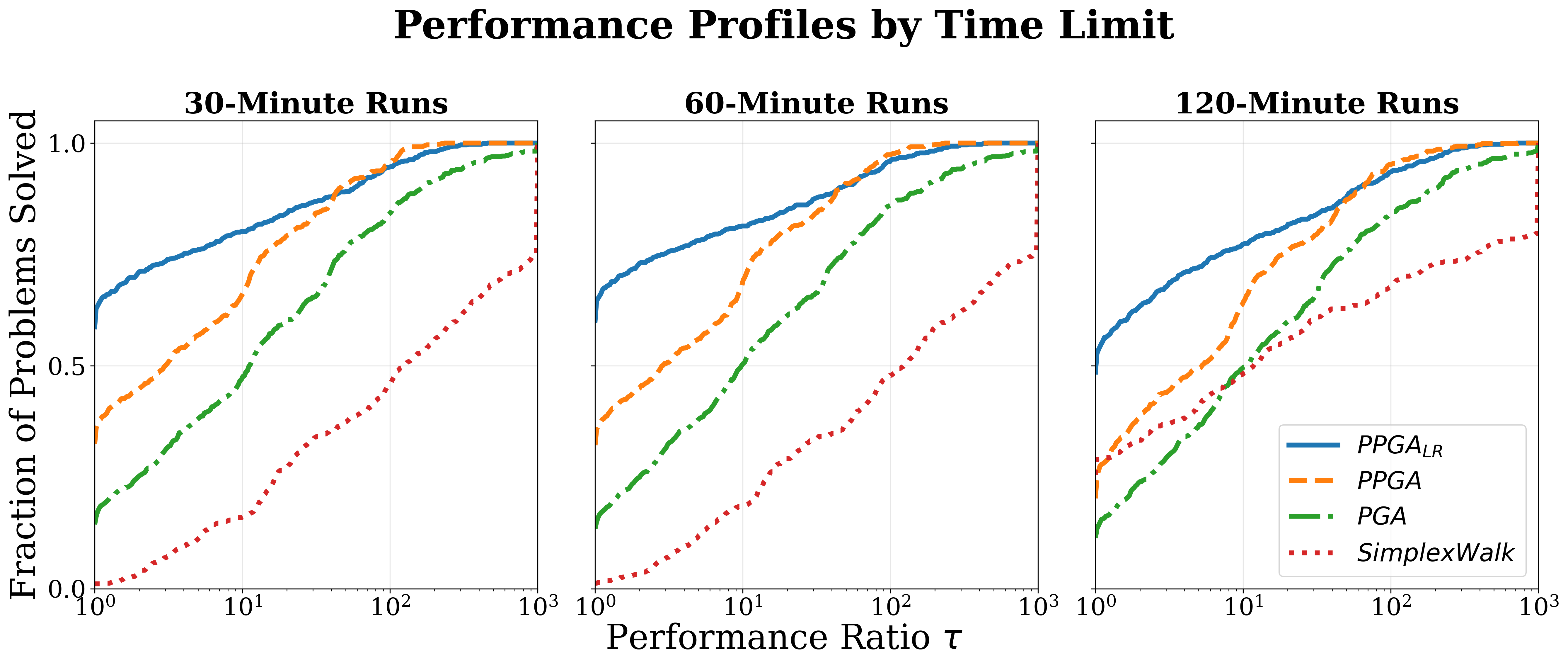}
    \caption{Comparison of algorithm performance over all instances by varying time limit.}
    \label{fig:overallbytime}
\end{figure}

\begin{figure}[h!]
    \centering
    \includegraphics[width=1.0\linewidth]{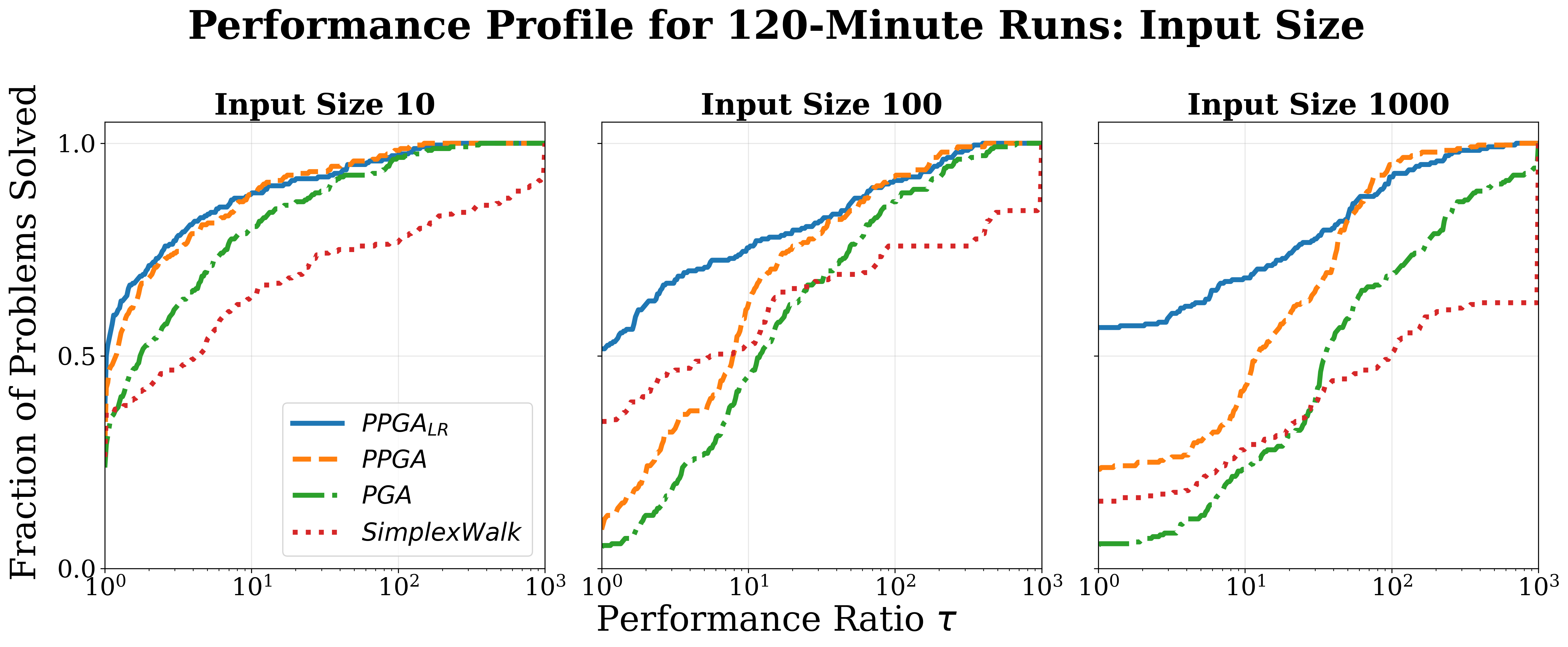}
    \caption{Comparison of algorithm performance by varying input size in 120-minute runs.}
    \label{fig:input}
\end{figure}

Figure \ref{fig:overallbytime} shows the overall performance profiles across all setups after the methods have been run for 30, 60, and 120 minutes. 
The other two plots focuses on the results after 120 minutes. 
Figure \ref{fig:input} shows the performance profiles from partitioning the instances according to the input size of the neural network. 
Figure \ref{fig:depthbywidth} partitions the instances according to the depth and the width of the neural network.

\subsection{Analysis}

In general, the relationship between the three gradient algorithms is similar across different time budgets, with \ppga{} and \ppgalr{} outperforming the vanilla \pga{}. 
However, the gap between the first \ppga{} and vanilla \pga{} shrinks for the largest width (10000). Since width dominates the complexity of gradient steps, it is possible that \ppga{} algorithms do not gain much in performance advantage in the early steps compared to the vanilla algorithm. Therefore, the advantage that we observe is likely gained in the later steps, when the vanilla algorithm is captured in some local optimum while \ppga{} escapes through perturbation.

\begin{figure}[hb!] 
    \centering
    \includegraphics[width=1.0\linewidth]{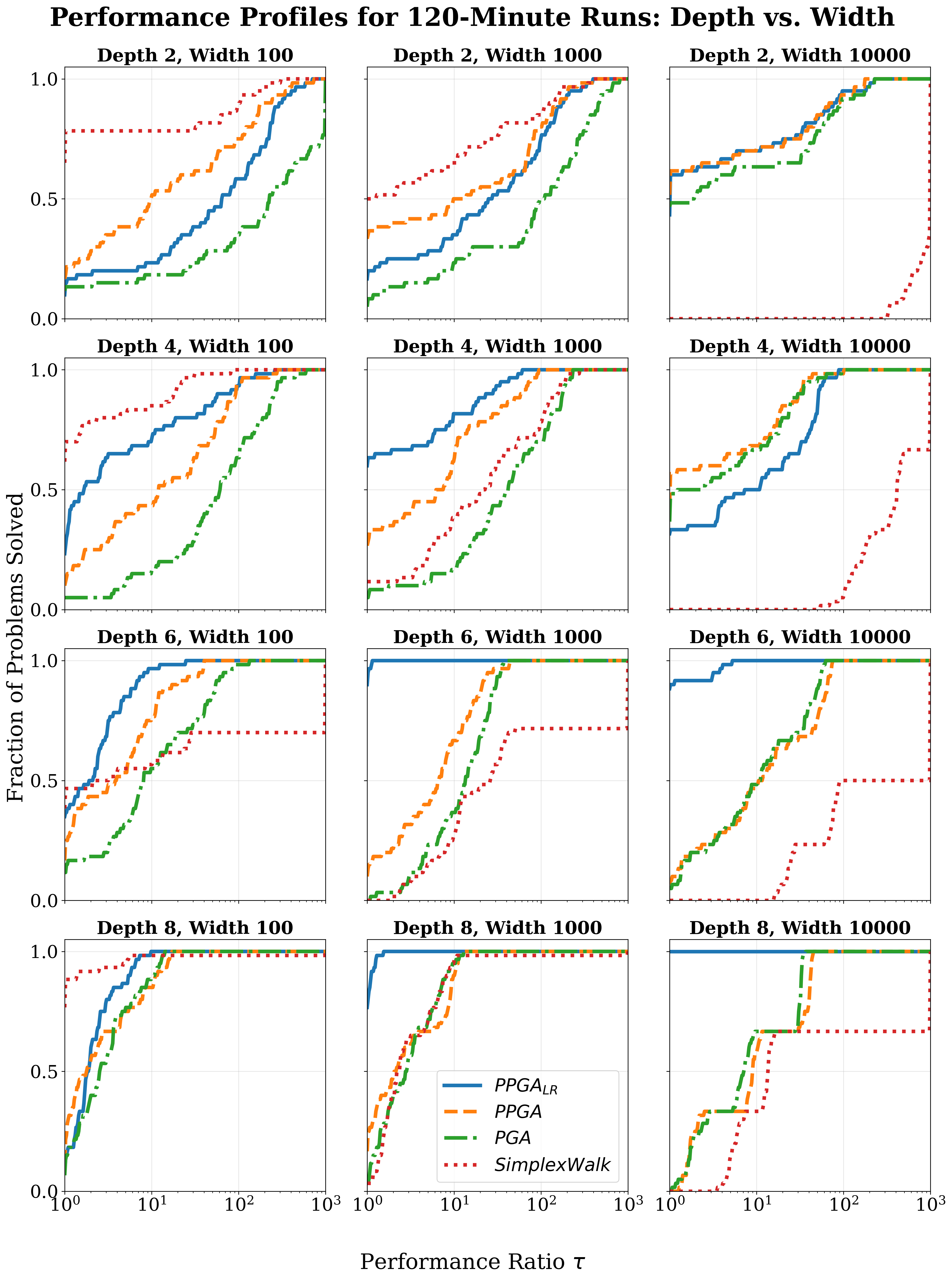}
    \caption{Comparison of performance by varying depth and width in 120-minute runs.}
    \label{fig:depthbywidth}
\end{figure}

While \ppgalr{} is overall better than \ppga{} in the aggregate of instances from Figure~\ref{fig:overallbytime}, 
the advantage of \ppgalr{} is due to the instances with larger dimensions:
\begin{itemize}
    \item From Figure~\ref{fig:input}, we see that \ppgalr{} matches \ppga{} for the smallest input (10) but does significantly better than \ppga{} for larger input sizes (100 and 1000).
    \item From the rows of Figure~\ref{fig:depthbywidth}, \ppgalr{} is generally better for larger widths.
    \item From the columns of Figure~\ref{fig:depthbywidth}, we see that \ppgalr{} is worse for the smallest depth (2), but that it dominates the results for deeper networks (6 and 8).
\end{itemize}

Hence, larger input size and depth seem to be the most determinant factors for \ppgalr{} performing better than \ppga{}.
Larger width also contributes in some cases. 
Those conditions conform with the cases in which neural networks tend to have more linear regions, 
which has inspired the design of \ppgalr{} in the first place. 
Moreover, 
as indicated by Equation~\ref{eq:grad} and given the weight distribution, 
neural networks with greater depth are bound to have smaller gradients. 

Meanwhile, notably, \sw{} gets comparatively better with longer runtimes, as shown in Figure~\ref{fig:overallbytime}. 
That conforms with the intuition that it takes longer to converge due to the more costly steps, but that each step tends to provide greater improvements. 
Indeed, 
we observe \sw{} dominating in the three scenarios with smallest depth and width in Figure~\ref{fig:depthbywidth}.


\subsection{A Case Study on Adaptive Stepsize for \ppgalr{}}

\begin{figure}[b]
    \centering
    \includegraphics[width=1.0\linewidth]{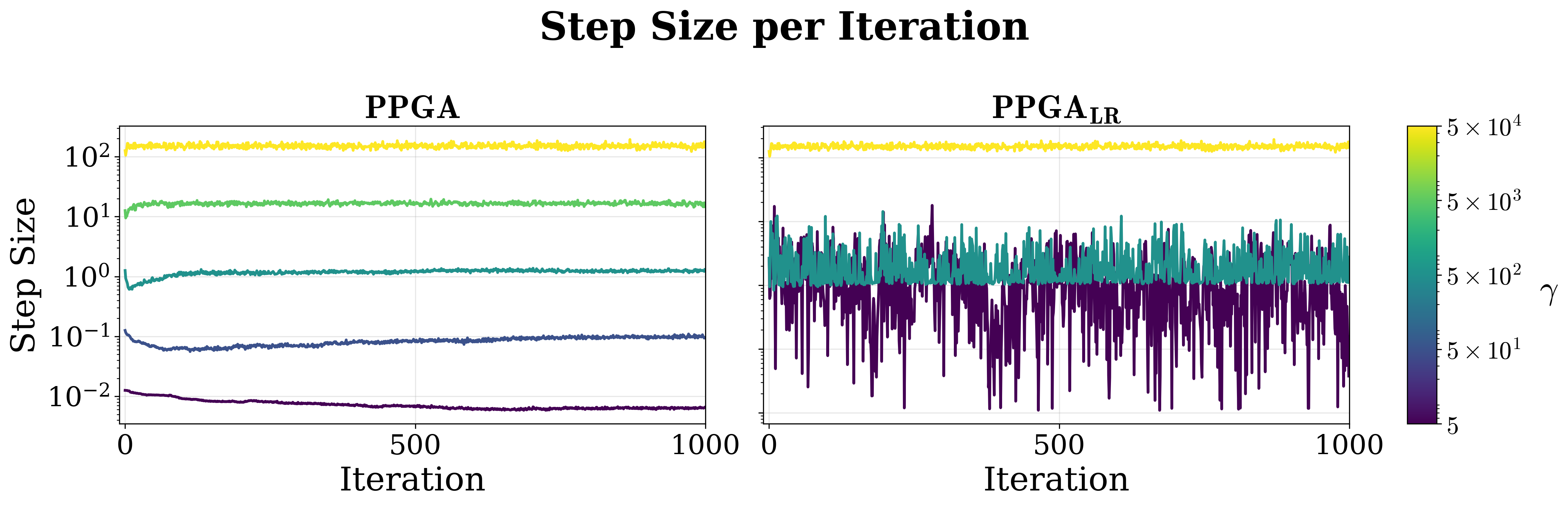}
    \caption{Step size in log scale over iterations for \ppga{} (Left) and \ppgalr{} (Right) for each learning rate.}
    \label{fig:stepsize}
\end{figure}

\begin{figure}[t]
    \centering
    \includegraphics[width=1.0\linewidth]{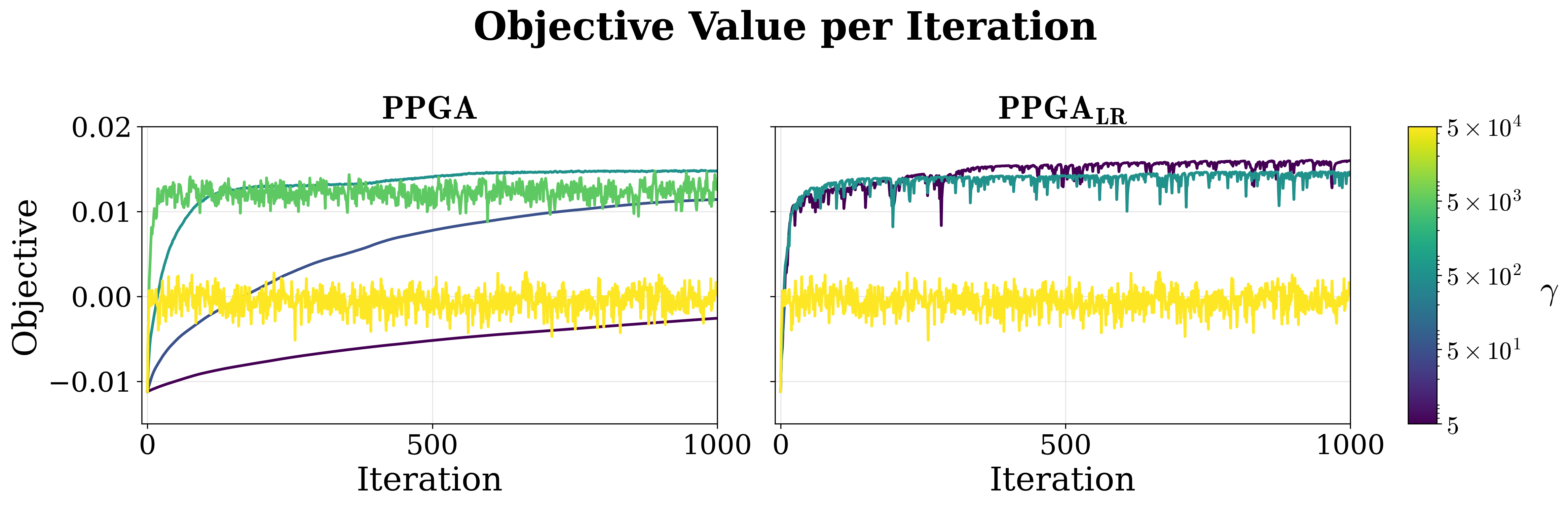}
    \caption{Solution value over iterations for \ppga{} (Left) and \ppgalr{} (Right) for each learning rate.}
    \label{fig:objective}
\end{figure}

To better understand why \ppgalr{} dominates in deeper instances, 
we present a case study on a randomly generated instance with $(n_0 = 1000, d=6, w=1000)$ and seed $s = 30$. We ran $1000$ iterations of \ppga{} with five learning rates, $\gamma\in\{5, 50, 500, 5000, 50000\}$.  We also ran $1000$ iterations of \ppgalr{}, but with only three different learning rates, $\gamma\in\{5, 500, 50000\}$. We use fewer learning rates with \ppgalr{} because there is a greater overlap in step sizes and objective values across learning rates in this case, making it harder to visually distinguish them. Figure \ref{fig:stepsize} shows the influence of the learning rate $\gamma$ on the step size over time. 
Figure \ref{fig:objective} shows the influence of the learning rate $\gamma$ on the solution value.

Our three main observations are the following:
\begin{itemize}
    \item In \ppga{}, when we increase the learning rate $\gamma$, the step size changes proportionally.  In \ppgalr{}, on the other hand, the step size is less affected by the learning rate and only the size of the smallest steps are clipped. 
    \item In \ppga{}, the range of step sizes for a specific learning rate is small. In \ppgalr{}, on the other hand, the valve mechanism defined by Equation~(\ref{eq: valve}) allows for a larger range of step sizes regardless of learning rate. That seems to lead to similar convergence for step sizes in different orders of magnitude.
    \item Since there is no mechanism designed for avoiding excessively large step sizes, both algorithms perform similarly worse when $\gamma = 50000$. We could potentially address this by designing another valve and truncating large steps. 
\end{itemize}
Hence, \ppga{} is more sensitive to the learning rate and does not converge well with lower learning rates ($\gamma \in \{{5}, {50}\}$). In turn, \ppgalr{} is more robust and converges well with lower learning rates in deep neural networks, as observed in Figure \ref{fig:depthbywidth}.

\section{Conclusion} \label{sec:conclusion}

In this paper, we proposed Gradient Walk: a new approach for local search based on gradient steps to optimize over the piecewise affine landscape of a trained ReLU network function within a polytope. 
We also proposed two algorithms for this approach, \ppga{} and \ppgalr{}. 
\ppga{} is based on projected gradient steps and empirically-crafted perturbations. 
\ppga{} is designed to search for a feasible solution while avoiding entrapment in arbitrary local optima constructed by the special landscape of ReLU networks. 
\ppgalr{} extends \ppga{} by using 
adaptive gradient steps that stretch further when the step size is too small relative to the local affine landscape.  
\ppgalr{} is designed to circumvent the 
issue of diminishing step sizes caused by the combination of a misaligned learning rate with small gradients observed in deep ReLU networks.  
We report that \ppgalr{} has a better performance than \ppga{} in neural networks with larger dimensions, whereas both outperform the baseline projected gradient ascent algorithm \pga{} and compare favorably against the LP Walk algorithm \sw{}, 
which has previously shown better results than other alternatives in the benchmark setting~\cite{tong2024optimization}.

In a nutshell, 
we double down on the trend of using specialized local search with lower per-iteration cost to find better solutions in constraint learning. 
We achieve that with gradient steps leveraging the local structure in ReLU networks.

Given the emergence of constraint learning frameworks~\cite{bergman2022janos,lueg2021relumip,ceccon2022omlt,fajemisin2023ocl,maragno2023mixed,gurobi2025ml,turner2025pyscipoptml} 
for applications in operations research, 
the development of new algorithms to timely find better solutions is of particular importance.

In future work, 
we intend to devise a variant of \ppga{} that is also robust against large learning rates, 
explore how other characteristics of the landscape of ReLU networks can be leveraged to produce better algorithms, 
apply the Gradient Walk approach to solve other types of constraint learning models, 
and explore its performance on specific families of functions approximated by neural networks.

\begin{credits}

\subsubsection{\discintname}
The authors have no competing interests to declare that are
relevant to the content of this article. 
\end{credits}


%
%
%
\FloatBarrier
\bibliographystyle{splncs04}
\bibliography{cites}

\end{document}